\documentclass{amsart}

\usepackage{amssymb,amsfonts,amsmath}

\usepackage{amssymb,amsbsy,amsthm,amsmath,graphicx,epsfig}

\newtheorem{theorem}[subsection]{Theorem}
\newtheorem{proposition}[subsection]{Proposition}
\newtheorem{lemma}[subsection]{Lemma}
\newtheorem{corollary}[subsection]{Corollary}

\newtheorem{conjecture}[subsection]{Conjecture}

\newtheorem{example}[subsection]{Example}

\newtheorem{remark}[subsection]{Remark}

\numberwithin{equation}{section} \theoremstyle{plain}

\renewcommand{\leq}{\leqslant}
\renewcommand{\geq}{\geqslant}
\newsavebox{\proofbox}
\savebox{\proofbox}{\begin{picture}(7,7)  \put(0,0){\framebox(7,7){}}\end{picture}}

\newcommand\Z{\mathbb{Z}}
\newcommand\R{\mathbb{R}}

\newcommand\N{\mathbb{N}}

\def\g{\mathfrak{g}}

\newcommand\eps{\varepsilon}
\newcommand\id{\operatorname{id}}

\begin{document}

\title[Rate of convergence, asymptotic cones and subFinsler geometry]{On the rate of convergence to the asymptotic cone for nilpotent groups and subFinsler geometry}

\author{Emmanuel Breuillard}
\address{Laboratoire de Math\'ematiques\\
B\^atiment 425, Universit\'e Paris Sud 11\\
91405 Orsay\\
FRANCE}
\email{emmanuel.breuillard@math.u-psud.fr}

\author{Enrico Le Donne}
\address{Laboratoire de Math\'ematiques\\
B\^atiment 425, Universit\'e Paris Sud 11\\
91405 Orsay\\
FRANCE}
\email{ledonne@msri.org, enrico.ledonne@math.u-psud.fr}

\keywords{nilpotent groups, subRiemannian geometry, polynomial growth, asymptotic cone,
Gromov-Hausdorff convergence, abnormal geodesics}
\subjclass{53C17, 22E25}



\begin{abstract} Addressing a question of Gromov, we give a rate in Pansu's theorem about the convergence in Gromov-Hausdorff metric of a finitely generated nilpotent group equipped with a left-invariant word metric scaled by a factor $\frac{1}{n}$ towards its asymptotic cone.  We show that due to the possible presence of abnormal geodesics in the asymptotic cone, this rate cannot be better than $n^{\frac{1}{2}}$ for general non-abelian nilpotent groups. As a corollary we also get an error term of the form $vol(B(n))=cn^d + O(n^{d-\frac{2}{3r}})$ for the volume of Cayley balls of a nilpotent group with nilpotency class $r$. We also state a number of related conjectural statements. \end{abstract}

\maketitle

\tableofcontents

\section{Introduction}

In his fundamental paper on groups with polynomial growth \cite{Gromov-polygrowth}, Gromov observed that Cayley graphs of finitely generated groups with polynomial growth, when viewed from afar, admit limits that are Lie groups endowed with a certain left-invariant geodesic metric. This was a simple but basic step in his proof that groups with polynomial growth admit nilpotent subgroups of finite index.\\

In his thesis \cite{pansu}, Pansu established that, if we start with the Cayley graph of a nilpotent group, then there is a unique limit. In other words, the sequence of metric spaces $\{Cay(\Gamma,\frac{1}{n}\rho_S)\}_{n \in \N}$ of scaled down Cayley graphs of the nilpotent group $\Gamma$ with generating set $S$ and word metric $\rho_S$ converges in the Gromov-Hausdorff topology \cite{Gromov-Metric-structures} towards a certain explicit left-invariant subFinsler metric on a nilpotent Lie group: the \emph{Pansu limit metric} on the asymptotic cone of $\Gamma$ (see Figure 1).\\

The goal of this note is to study the rate of convergence in Pansu's theorem and give quantitative estimates. This question was posed by Gromov in \cite[\S 2C Remark 2C2(a)]{gromov-asymp}. It requires approximating (with explicit bounds) word metric geodesics in $\Gamma$ with subFinsler geodesics in the Lie group that is the asymptotic cone of $\Gamma$ and vice-versa. This problem is intimately connected to the underlying geometry of nilpotent Lie groups endowed with left-invariant subFinsler metrics.\\

One of our key findings is that the quality of the error term in Pansu's theorem is related to the presence or the absence of so-called \emph{abnormal geodesics} in the asymptotic cone. These geodesics do not exist in classical Riemannian or Finsler geometry, but are typical in subRiemannian or subFinsler geometry  (see \cite{Montgomery}). We show that their presence worsens the error term in the convergence to the asymptotic cone for general nilpotent groups. For example, if $\Gamma=\Z^n$, the free abelian group of rank $n$, or if $\Gamma=H_{2n+1}(\Z)$, the $2n+1$-dimensional Heisenberg group, equipped with a word metric, then the asymptotic cone of $\Gamma$ bears no abnormal geodesics, and it can be shown that the convergence to the asymptotic cone is best possible, namely with a rate $\frac{1}{n}$. On the other hand, if $\Gamma=H_3(\Z) \times \Z$, the speed of convergence towards the asymptotic cone depends on the word metric, and, while for some generating sets the speed may be optimal, i.e., in $\frac{1}{n}$, we show that, for some choices of generating sets, the rate of convergence is no faster that $n^{-\frac{1}{2}}$. This is due to the fact that the asymptotic cone $H_3(\R) \times \R$ admits abnormal geodesics in the direction of the second $\R$ factor (see Figure 2 and Figure 3). In particular, we exhibit two word metrics on $H_3(\Z) \times \Z$ that have isometric asymptotic cones and yet are not $(1,C)$-quasi-isometric for any constant $C>0$, answering negatively a related question raised by Burago and Margulis \cite{Oberwolfach2006}.\\

For general nilpotent groups, we obtain a rate of convergence in $O(n^{-\frac{2}{3r}})$ if $r>2$, and $O(n^{-\frac{1}{2}})$ if $r=2$, where $r$ is the nilpotency class of $\Gamma$. The case $r=2$ is a simple consequence of a result of Stoll \cite{stoll} and we use it to obtain the $\frac{2}{3r}$ exponent for general $r$. The above example on $H_3(\Z) \times \Z$ shows that this  error term is sharp for groups of nilpotency class $2$. However, it is likely not to be sharp for groups of nilpotency class $3$ or more and we conjecture that the exponent $\frac{1}{2}$ holds for all nilpotent groups and is therefore independent of $r$.\\

The rate of convergence in Gromov-Hausdorff metric towards the asymptotic cone implies a rate of convergence in the volume asymptotics of Cayley balls. In particular, we obtain as a corollary that, for every nilpotent group $\Gamma$ with generating set $S$ ($S$ finite, $S=S^{-1}$ and $1\in S$), we have
that 
\begin{equation}\label{asym}
|S^n|=c_S n^d +O_S(n^{d-\beta_r}),\quad \text{as } n\to +\infty,
\end{equation}
 where $\beta_r=\frac{2}{3r}$ for $r>2$ and $S^n$ is the ball of radius $n$ in the word metric and $c_S>0$ a constant.\\

Stoll had showed in \cite{stoll} that one can take $\beta_2=1$ for groups of nilpotency class at most $2$, but nothing seemed known for higher-step groups, even though it is a folklore conjecture that $\beta_r=1$ should hold for all $r$. We also note that the very fact that an error term exists at all is also a distinctive feature of nilpotent groups. 
Indeed, there is a class of (non-discrete) groups of polynomial growth (of the form $\R^2 \rtimes_\theta \Z$, where $\Z$ acts by an irrational rotation with angle $\theta$), which are solvable but not nilpotent, for which it can be shown that the volume asymptotics (here $vol(B(n)) \simeq cn^3$) admit no error term whatsoever if $\theta$ is chosen carefully \cite{breuillard}.\\ 


In this note we present the aforementioned results and give some information on their proofs, although full details will be given elsewhere due to lack of space.\\

The note is organized as follows. First, we state our main result about the rate of convergence to the asymptotic cone, Theorem \ref{main-theorem} below, and explain the strategy of the proof and its main ingredients. In the second part of the paper, we present the example showing the sharpness of the exponent $\frac{1}{2}$ for step-$2$ groups (Theorem \ref{noquasi} and Proposition \ref{quasi}). We then  describe some ideas in the proof of the sharpness result, in particular a detailed study of the geometry of subFinsler metrics on the Heisenberg group and its direct product with $\R$. Finally, in the last part of the note, we prove the volume estimate $(\ref{asym})$, discuss the volume asymptotics conjecture, and its relation to some well-known conjectures in subRiemannian geometry.

\begin{figure}\label{heisenball}
\begin{center}
\includegraphics[scale=.4]{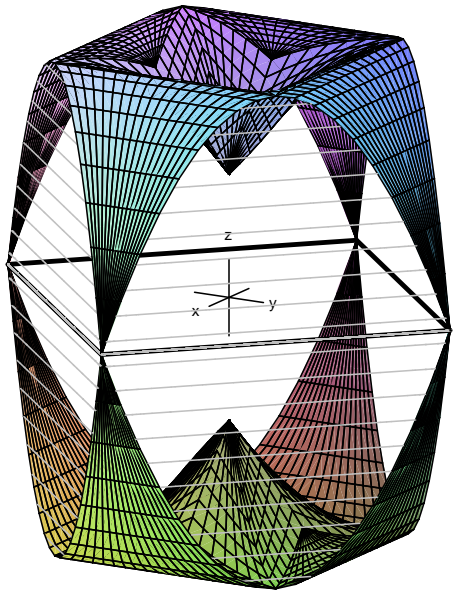}
\caption{
 The unit ball for the Pansu limit metric $d_3$ of the Cayley graph of the discrete Heisenberg group $H_3(\Z)$ with standard generators (see later in Section \ref{fine} for an explicit formula for $d_3$).}
\end{center}
\end{figure}

\section{The rate of convergence to the asymptotic cone}
We now recall Pansu's description of the asymptotic cone of a nilpotent group and state our main theorem. Let $\Gamma$ be a torsion-free nilpotent group generated by a finite set $S$ (we assume $1 \in S$ and $S^{-1}=S$). It is well-known \cite{Raghunathan} that $\Gamma$ embeds as a co-compact discrete subgroup of a simply connected nilpotent Lie group $G$, its Malcev closure. Let $\g$ be the Lie algebra of $G$. One can associate to $\g$ another nilpotent Lie algebra, called the \emph{graded Lie algebra of $\g$} and denoted by $\g_\infty$ by setting $$\g_\infty=\oplus_{i \geq 1} \g^{(i)}/\g^{(i+1)},$$
where $\g^{(i+1)}=[\g,\g^{(i)}]$ is the central descending series of $\g$ and where the Lie bracket is defined in the natural way by the formula $[x\textnormal{ mod }\g^{(i+1)},y\textnormal{ mod }\g^{(j+1)}]=[x,y]\textnormal{ mod }\g^{(i+j+1)}$ for $x \in \g^{(i)}$ and $y \in \g^{(j)}$.\\

The exponential map establishes a diffeomorphism between $\g$ and the Lie group $G$ and between $\g_\infty$ and a simply connected Lie group $G_\infty$. We denote the group law on $G_\infty$ by $x*y$ to distinguish it from the group law on $G$, for which we simply write $x\cdot y$. The projection Lie algebra homomorphism
\begin{equation}\label{proj}\pi: \g \to \g/[\g,\g],
 \end{equation}
 lifts to a group homomorphism $\pi:G \to \g/[\g,\g]$ which, by a slight abuse of notation, we also denote by $\pi$. The image $\pi(\Gamma)$ is a discrete co-compact subgroup of the vector space $\g/[\g,\g]$ generated by $\pi(S)$. In particular $\pi(S)$ defines a norm on $\g/[\g,\g]$ whose unit ball is the convex hull of $\pi(S)$. We call this norm $\|\cdot\|_{S}$ the \emph{Pansu limit norm} of the pair $(\Gamma,S)$.\\

Viewing $\g/[\g,\g]$ as a subspace of the graded Lie algebra $\g_\infty$, we may define a left-invariant subFinsler metric on $G_\infty$ with horizontal subspace $\g/[\g,\g]$ and norm $\|\cdot\|_{S}$. In other words, there is a left-invariant geodesic metric $d_\infty$ on $G_\infty$ defined by
$$d_\infty(x,y):=\inf\{ L(\gamma)\},$$
where the infimum is taken over all piecewise linear horizontal paths $\gamma$ from $\gamma(0)=x$ to $\gamma(1)=y$, and $L(\gamma)$ is the length of $\gamma$ measured using the the Pansu limit norm $\|\cdot \|_S$. A piecewise horizontal path is by definition the concatenation of finitely many segments of one parameter subgroups of $G_\infty$ of the form $\{\exp(tX)\}_{t\in [0,T]}$ for some $X \in \g/[\g,\g]$.\\

The distance $d_\infty$ is geodesic and left-invariant by definition. We call it subFinsler, because the norm $\|\cdot\|_S$ is not a Euclidean norm but a polyhedral norm. It is not a Finsler metric however (if $\g$ is non-abelian), because the norm is only defined on a subspace of the Lie algebra. It can be checked that this subspace $\g/[\g,\g]$ generates the whole Lie algebra. This implies (Chow's theorem \cite{Montgomery}) that every two points in $G_\infty$ can be joined by a piecewise linear horizontal path and thus $d_\infty$ is well-defined. We note finally that it can be shown \cite{Berestovskii} that any left-invariant geodesic metric on a connected Lie group is a subFinsler metric for some norm on a generating subspace of the Lie algebra.\\

We call the distance $d_\infty$ the \emph{Pansu limit metric} of $(\Gamma,S)$. Pansu showed in his thesis that $(G_\infty,d_\infty)$ is the asymptotic cone of $(\Gamma,\rho_S)$. More precisely he showed that
the sequence of renormalized Cayley graphs $Cay(\Gamma,\frac{1}{n}\rho_S)$ of $\Gamma$ converges towards $(G_\infty,d_\infty)$ in the Gromov-Hausdorff topology on pointed metric spaces (based at $\id$). Here $\rho_S$ is the left-invariant word metric on $\Gamma$ defined by $\rho_S(x,y):=\inf\{n \in \N; x^{-1}y \in S^n\}$.\\

For any $R>0$, set $X_\infty(r):=(B_{d_\infty}(\id,R),d_\infty)$, where $B_{d_\infty}(\id,R)$ is the closed ball of radius $R$ in $(G_\infty,d_\infty)$. Similarly, we set $X_n(R):=(B_{S}(\id,Rn),\frac{1}{n}\rho_S)$, where $B_{S}(\id,Rn)$ is the closed ball of radius $Rn$ in the Cayley graph $(\Gamma,\rho_S)$. Then\\

\begin{theorem}[Pansu \cite{pansu}] For every $R>0$, $X_n(R)$ converges to $X_\infty(R)$ in the Gromov-Hausdorff topology.\\
\end{theorem}

Recall that the Gromov-Hausdorff metric between any two compact metric spaces $(X,d_X)$ and $(Y,d_Y)$ is defined as
$$d_{GH}(X,Y):=\inf\{d_{\rm Haus,Z}(X,Y); Z=X \sqcup Y, d_Z\textnormal{ admissible}\},$$
where an admissible metric $d_Z$ on the disjoint union $Z=X \sqcup Y$ is a distance on $Z$ which coincide with $d_X$ on $X$ and with $d_Y$ on $Y$.\\

Here is our main result:
\begin{theorem}\label{main-theorem} There is a positive constant $\alpha_r>0$ depending only on the nilpotency class $r$ of $\Gamma$ such that. as $n \to +\infty,$
$$d_{GH}(X_n(1),X_\infty(1)) = O_S(n^{-\alpha_r}).$$
Moreover one can take $\alpha_1=1$, $\alpha_2=\frac{1}{2}$ and $\alpha_r=\frac{2}{3r}$ if $r>2$.\\
\end{theorem}

Note that $X_n(1)$ as a set is the ball of radius $n$ for the word metric $\rho_S$ on $\Gamma$. Note also that by scaling, this also implies that $$d_{GH}(X_n(R),X_\infty(R))=O_S(R^{1-\alpha_r}n^{-\alpha_r}),$$ for every $R>0$.\\

Theorem \ref{main-theorem} addresses a question of Gromov \cite[Remark 2C2(a)]{gromov-asymp}. Speaking of the Gromov-Hausdorff distance between $X_n(1)$ and $X_\infty(1)$ he says ``any bound on these distances would be a pleasure to have in our possession, even in the case of word metrics on $\Gamma$''. \\

Although we have proved Theorem \ref{main-theorem} for word metrics on $\Gamma$ only, using Burago's theorem \cite{burago}, one can adapt our arguments and prove a similar result (at least with $\alpha_r\geq 1/2r$) for $\Gamma$-invariant Riemannian metrics on the Malcev closure $G$ of $\Gamma$ and more generally for $\Gamma$-invariant coarsely geodesic metrics on $G$. We will not pursue this here.\\

Our result is sharp for $r=1,2$. The sharpness in case $r=2$ is discussed below and is related to a conjecture of Burago and Margulis \cite{Oberwolfach2006} and to the presence of abnormal geodesics in the asymptotic cone. However, we believe that the exponent is not sharp for $r>2$ and that the exponent $\frac{1}{2}$ holds for all $r\geq 2$. Proving this would require a deeper understanding of the subFinsler geodesics of $d_\infty$ that we have so far. The proof of Theorem \ref{main-theorem} is to a large extent an effectivization of Pansu's proof, where we have to replace several compactness arguments used by Pansu by effective arguments with explicit bounds. The exponent $\frac{1}{2}$ when $r=2$ is deduced from a theorem of Stoll \cite{stoll} and Stoll's result is also used in order to obtain the exponent $\frac{2}{3r}$ for $r>2$.\\

\section{Comparison of metrics on a nilpotent Lie group}
In order to prove Theorem \ref{main-theorem}, one needs to first have some understanding of the large scale behavior of subFinsler metrics on nilpotent Lie groups. Given two left-invariant subFinsler metrics $d_1$ and $d_2$ on a simply connected nilpotent Lie group $G$, the following gives a simple criterion for when they are asymptotic and also gives an error estimate.\\

\begin{proposition}\label{subFinsler-comp} Given two left-invariant subFinsler metrics $d_1$ and $d_2$ on $G$, the following are equivalent:
\begin{enumerate}
\item $\frac{d_1(\id,g)}{d_2(\id,g)} \rightarrow 1$, as $g \to \infty$ in $G$,
\item the projection under $\pi$ (see $(\ref{proj})$) of the unit balls of $\|\cdot\|_1$ and $\|\cdot\|_2$ coincide,
\item $|d_1(\id,g) - d_2(\id,g)|=O(d_1(\id,g)^{1-\frac{1}{r}})$, as $g \to \infty$ in $G$.
\end{enumerate}
\end{proposition}

The norms $\|\cdot\|_1$ and $\|\cdot\|_2$ are the norms used to define $d_1$ and $d_2$. respectively, as recalled in the preceding section. Similarly, one can prove that $(G,d_1)$ and $(G,d_2)$ have isometric asymptotic cones if and only if the normed vector space $\g/[\g,\g]$ endowed with the projection of $\|\cdot \|_1$ is isometric to the same space endowed with the projection of $\|\cdot\|_2$.\\

Item $(ii)$ in the above proposition can also be replaced with `the projection under $\pi$ of the unit balls of $d_1$ and $d_2$ coincide'. In this form the proposition applies also to word metrics, instead of subFinsler metrics, that is the pseudo-distances on $G$ of the form $d_U(x,y)=\inf\{n \in \N; x^{-1}y \in U^n\}$, where $U$ is a bounded symmetric neighborhood of $\id$ in $G$.\\

The above dealt only with $G$-invariant metrics. For Theorem \ref{main-theorem}, one also needs to consider metrics on the discrete co-compact subgroup $\Gamma$, or more generally $\Gamma$-invariant metrics on $G$. This is more involved, because small errors made in discretizing at various places in the group can accumulate and generate a large error in the end.\\

Let $G$ be a simply connected nilpotent Lie group and $\Gamma$ a discrete co-compact subgroup of $G$. Let $\rho_S$ be the left-invariant word metric on $\Gamma$ associated to a finite symmetric set $S$.  The following result is very closely related to Theorem \ref{main-theorem} and gives an estimate of the error term between $\rho_S$ and $G$-invariant subFinsler metrics on $G$.\\

\begin{theorem} \label{main-theorem-second} Let $d$ be a left-invariant subFinsler metric on $G$. The following are equivalent:
\begin{enumerate}
\item $\frac{d(\id,\gamma)}{\rho_S(\id,\gamma)} \rightarrow 1$, as $\gamma \to \infty$ in $\Gamma$,
\item the projection under $\pi$ (see $(\ref{proj})$) of the unit ball of $\|\cdot\|_d$ is the convex hull of $\pi(S)$,
\item $|d(\id,\gamma) - \rho_S(\id,\gamma)|=O(d(\id,\gamma)^{1-\alpha_r})$, as $\gamma \to \infty$ in $\Gamma$, where $\alpha_r$ is as in Theorem \ref{main-theorem}.
\end{enumerate}
\end{theorem}

This is consistent with Proposition \ref{subFinsler-comp} because $\alpha_r \geq \frac{1}{r}$. One may wonder however if, given $S$ and $\rho_S$, there is a distinguished $G$-invariant subFinsler metric on $G$ that best approximates $\rho_S$. A good candidate for this is the following choice of subFinsler metric, which we call the \emph{Stoll metric} associated to $(\Gamma,S)$.\\

Identifying $G$ with its Lie algebra $\g$ via the exponential map, we may view the finite symmetric set $S$ as a subset of $\g$ and take its convex hull. It spans a vector subspace $V_S$ of $\g$. Let $\|\cdot\|_S$ be the norm on $V_S$ whose unit ball is the convex hull of $S$. Then $\|\cdot\|_S$ induces a left-invariant subFinsler metric $d_S$ on $G$, which we call the Stoll metric.\\

\begin{conjecture}\label{stoll-conj} There is a constant $C=C(S)>0$ such that
$$|\rho_S(\id,\gamma) - d_S(\id,\gamma)| \leq C,$$
for all $\gamma \in \Gamma$.
\end{conjecture}

In a very elegant short paper \cite{stoll} M. Stoll  proved this claim\footnote{Stoll's definition of $d_S$ is slightly different, but it is easy to see that it yields the same distance.} when the nilpotency class of $\Gamma$ is at most $2$. It remains an open problem for $r>2$. In fact, even the existence of some $G$-invariant metric on $G$ that lies at a bounded distance from $\rho_S$ seems unknown.\\

Of course, by Proposition \ref{subFinsler-comp}, the Stoll metric and the $G$-invariant Pansu limit metric (induced by $\|\cdot\|_\infty$) are asymptotic, but the Stoll metric seems to capture the finer behavior of the word metric.\\

The fact that Conjecture \ref{stoll-conj} holds for $r=2$ has the following simple consequence. Let $\pi_2 : G \to G/G^{(3)}$ be the projection homomorphism modulo $G^{(3)}=[G,[G,G]]$ and let $d_{\pi_2(S)}$ be the Stoll metric on $G/G^{(3)}$ associated to $\pi_2(S)$.\\

\begin{lemma} \label{stoll-lemma}There is a constant $C=C(S)>0$ such that for every $u \in G/G^{(3)}$, there is $\gamma \in \Gamma$ with $\rho_S(\id,\gamma) \leq d_{\pi_2(S)}(\id,u)$ and such that $d_{\pi_2(S)}(\pi_2(\gamma), u) \leq C$.\\
\end{lemma}

This lemma will be helpful in the proof of Theorem \ref{main-theorem} in order to approximate a subFinsler geodesic by a word metric geodesic.\\

\section{From continuous geodesics to discrete ones and back}

We now give a brief sketch of the proof of Theorem \ref{main-theorem}. For the sake of simplicity we will assume that the Lie algebra $\g$ is isomorphic to its graded $\g_\infty$. Additional technicalities arise if this is not the case, but they do not affect our convergence rates. In what follows we identify the Lie group with its Lie algebra via the exponential map. On $\g_\infty$ there is a natural one-parameter group $\{\delta_t\}_{t>0}$ of automorphisms called \emph{dilations}, and defined by the formula $\delta_t(x)=t^ix$ if $x \in \g^{(i)}/\g^{(i+1)}$. The dilations scale the subFinsler metric $d_\infty$ nicely and we have:
\begin{equation}\label{scaling}
d_\infty(\delta_t(x),\delta_t(y))=td_\infty(x,y).
\end{equation}

Theorem \ref{main-theorem} can be easily deduced from the following proposition. We recall that $B_S(\id,n)$ is the ball of radius $n$ in $\Gamma$ for the word metric $\rho_S$ and $d_\infty$ is the associated Pansu limit metric on the asymptotic cone of $\Gamma$.\\
\begin{proposition}\label{thmbis} Let $n \in \N$, then, for $\alpha_r$ as in Theorem \ref{main-theorem},
\begin{enumerate}
\item for every $\gamma \in B_S(\id,n)$, there is $y \in B_{d_\infty}(\id,1)$ such that $d_\infty(y,\delta_{\frac{1}{n}}(\gamma))=O(n^{-\frac{1}{r}})$, as $n\to \infty$,
 \item for every $x\in B_{d_\infty}(\id,1)$ there is $\gamma_x \in B_S(\id,n)$ such that $d_\infty(x,\delta_{\frac{1}{n}}(\gamma_x))=O(n^{-\alpha_r})$, as $n\to \infty$.
 \end{enumerate}
\end{proposition}

The two parts of the proposition are fairly distinct. The first one follows immediately from Proposition \ref{subFinsler-comp}, because $\rho_S(\id,\gamma) \leq d_S(\id, \gamma)$ for all $\gamma \in \Gamma$, where $d_S$ is the Stoll metric defined in the previous section. The second part of the proposition lies deeper as we need to approximate a $d_\infty$-geodesic in $G$ with a $\rho_S$-geodesic in $\Gamma$.\\

This is done by first splitting a $d_\infty$-geodesic in $G$ between $\id$ and $\delta_n(x)$ into $m \simeq n^{\frac{1}{3}}$ intervals of equal length, so that $\delta_n(x)=y_1 \cdot \ldots\cdot y_m$. Then one shows using the Campbell-Baker-Hausdorff formula combined with Lemma \ref{gronwall} below that $y_1 \cdot \ldots\cdot y_m$ can be approximated by $\pi_2(y_1)\cdot \ldots \cdot \pi_2(y_m)$ with only an error of order at most $\frac{n}{m^{2/r}}$. Here $\pi_2(y_i) \in G/G^{(3)}$ is viewed as an element of $G$ by viewing $\g/\g^{(2)} \oplus \g^{(2)}/\g^{(3)}$ as a subspace of $\g\simeq \g_\infty$. Finally one applies Lemma \ref{stoll-lemma}, which itself relies on Stoll's result \cite{stoll}, in order to approximate each $\pi_2(y_i)$ by a suitable element of $\Gamma$ and this ends the proof.\\

The approximation of $y_1\cdot \ldots\cdot y_m$ by $\pi_2(y_1)\cdot \ldots \cdot \pi_2(y_m)$, as well as the proof of Proposition \ref{subFinsler-comp} relies mainly on the Campbell-Baker-Hausdorff formula and the following simple fact, also used in \cite{pansu}, which is itself a version of the classical Gronwall's lemma. We note that a similar argument was used by the second-named author in \cite[Appendix]{LeDonne1}.\\

\begin{lemma}[Gronwall-type lemma]\label{gronwall}  Let $G$ be a Lie group. Let $\left\| \cdot \right\| $ be
some norm on the Lie algebra of $G$ and let $d_{e}(\cdot ,\cdot )$ be a
Riemannian metric on $G$. Then, for every $L>0$, there is a constant
$C=C(d_{e},\left\| \cdot \right\| ,L)>0$ with the following property. Assume
$\xi _{1},\xi _{2}:[0,1]\rightarrow G$ are two piecewise smooth paths in $G$ with $\xi _{1}(0)=\xi _{2}(0)=\id.$ Let $\xi _{i}^{\prime }\in
Lie(G)$ be the tangent vector pulled back at the identity by a left
translation of $G$. Assume that $\sup_{t\in [0,1]}\left\| \xi _{i}^{\prime
}(t)\right\| \leq L$, and that $\left\| \xi _{1}^{\prime
}(t)-\xi _{2}^{\prime }(t)\right\|_{ {L}^2([0,1])} \leq \varepsilon $. Then
\begin{equation*}
d_{e}(\xi _{1}(1),\xi _{2}(1))\leq C\varepsilon.
\end{equation*}
\end{lemma}

Note that the only reason for splitting the original $d_\infty$-geodesic into $m$ pieces is to be able to have long enough pieces so that the projection to $G/G^{(3)}$ can be well approximated by a $\rho_S$-geodesic using Lemma \ref{stoll-lemma}. Had we had Conjecture \ref{stoll-conj} at our disposal, both parts of Proposition \ref{thmbis} would follow directly from Proposition \ref{subFinsler-comp} and in the second part the coefficient $\alpha_r$ would be as good as in Proposition \ref{subFinsler-comp}, namely $\frac{1}{r}$.\\

Regarding Proposition \ref{subFinsler-comp} itself, we believe that the exponent $\frac{1}{r}$ can be replaced with $\frac{1}{2}$ even if $r>2$. We have checked this only when $r\leq 4$ so far. Combined with Conjecture \ref{stoll-conj}, having the exponent $1/2$
would imply that one could take $\alpha_r=\frac{1}{2}$ in Theorem \ref{main-theorem} also, at least for graded nilpotent groups.\\

Stoll's proof of Conjecture \ref{stoll-conj} for $r=2$ relies on a good understanding of geodesics for the $d_S$ metric, and in particular the fact that every point in $G$ can be joined by a $d_S$-geodesic that is piecewise horizontal linear with a bounded number of distinct linear pieces. Despite the fact that Lemma 3 in \cite{stoll} does not hold for $r\geq 3$, this latter property is likely to hold for all $r$.\\

\section{Abnormal geodesics, the Burago-Margulis conjecture and the sharpness of Theorem \ref{main-theorem}}
We now pass to the second part of this note, which concerns with the construction of a specific example showing that Theorem \ref{main-theorem} is sharp for groups of nilpotency class $2$. The fact that the best exponent $\alpha_2$ is $\frac{1}{2}$ instead of $1$ is somewhat surprising, not only because $\alpha_1=1$, but because of the fact that for the archetypal examples of step $2$ groups, namely the Heisenberg groups, the best exponent is also $1$, for any choice of generating set, see \cite{krat}.\\

In fact this issue is related to another surprising phenomenon of subRiemannian geometry, namely the existence of \emph{abnormal geodesics}, see \cite{Montgomery}. Loosely speaking, these are geodesics, say connecting two points $x$ and $y$, such that the endpoints of the $\eps$-variations of that geodesic do not cover a full ball (for a Riemannian metric in some chart) of radius $>C\eps$ around $y$, for some positive $C$. So typically, even when $x \neq y$, 
if the geodesic connecting $x$ and $y$ is abnormal, one will be able to find points $z$ near $y$ at distance say $\eps$ from $y$ in a Riemannian metric that are not at distance $d(x,y)+O(\eps)$ from $x$.\\

Abnormal geodesics do not exist in Riemannian geometry as one can see from the first variation formula, and they are a distinctive feature of subRiemannian geometry. A typical example is provided by segments of one parameter horizontal subgroups in the free nilpotent Lie group of step $2$ and rank at least $3$.\\

In this section, we consider the group $G=\R \times H_3(\R)$, the direct product of the $3$-dimensional Heisenberg group and the additive group of $\R$. We write an element of this group as $g=(v;x,y,z)$ if $g=(v,\exp(xX+yY+zZ))$, where $(X,Y,Z)$ forms a basis of the Lie algebra of $H_3(\R)$ such that $[X,Y]=Z$. It turns out that, if one puts a subFinsler product metric on $G$ starting with a subFinsler (but not Finsler) metric $d_3$ on $H_3(\R)$, say $d((0,\id),(t,g))=|t|+d_3(g,\id)$, then the segments $\{(t,\id)\}_{t \in [a,b]}$ are abnormal geodesics. Indeed, since $Z$ is the central direction in the Lie algebra of $H_3(\R)$, the points $(1,\exp(\eps Z))$ are at distance at least $1+c\sqrt{\eps}$ from $(0,\id)$ and not $1+O(\eps)$ as would be the case had $d_3$ been a Finsler (or Riemannian) metric on $H_3(\R)$.\\

For this reason, it is easy to cook up examples of subFinsler metrics $d_1$ and $d_2$ on $G$ that are asymptotic and yet have $|d_1(\id,g) - d_2(\id,g)| \gg d_1(\id,g)^{\frac{1}{2}}$ for arbitrarily large $g$, hence showing that the exponent $1-\frac{1}{r}$ in Proposition \ref{subFinsler-comp} is sharp when $r=2$. Indeed, simply take $\|\cdot\|_1=|v|+\max\{|x|,|y|,|z|\}$ and $\|\cdot\|_2=|v|+\max\{|x|,|y|\}$ for the norms defining the subFinsler metrics $d_1$ and $d_2$. Note that in this example $d_1$ is Finsler, while $d_2$ is only subFinsler with horizontal subspace $\R\times(\R X \oplus \R Y) $.\\

An analogous example in the discrete group $G(\Z):=H_3(\Z) \times \Z$ was given in \cite{breuillard} in order to disprove a conjecture of Burago and Margulis. Let $\rho_1$ and $\rho_2$ be the left-invariant word metrics on $G(\Z)$ induced by the finite generating sets
\begin{eqnarray*}
S_1:=\{(1;0,0,1)^{\pm 1} , (1;0,0,-1)^{\pm 1}, (0;1,0,0)^{\pm 1}, (0;0,1,0)^{\pm 1}\},
\\ S_2:=\{(1;0,0,0)^{\pm 1}, (0;1,0,0)^{\pm 1}, (0;0,1,0)^{\pm 1}\},
\end{eqnarray*}
respectively. Then it follows from Theorem \ref{main-theorem-second} that $\rho_1$ and $\rho_2$ are asymptotic in the sense that $\frac{\rho_1(\id,\gamma)}{\rho_2(\id,\gamma)}$ tends to $1$ as $\gamma$ tends to $+\infty$ and thus that $(G(\Z),\rho_1)$ and $(G(\Z),\rho_2)$ have isometric asymptotic cones. However, if $\gamma_n=(n;0,0,n)$, then one checks easily that $\rho_1(\id,\gamma_n)=n$ while $\rho_2(\id,\gamma_n)-n \gg \sqrt{n}$, see \cite{breuillard}. A picture of the unit ball of the asymptotic cone of $(G(\Z),\rho_1)$ is given in Figure 2.\\

  \begin{figure}\label{pallaA}
  \begin{center}
\includegraphics[scale=.45]{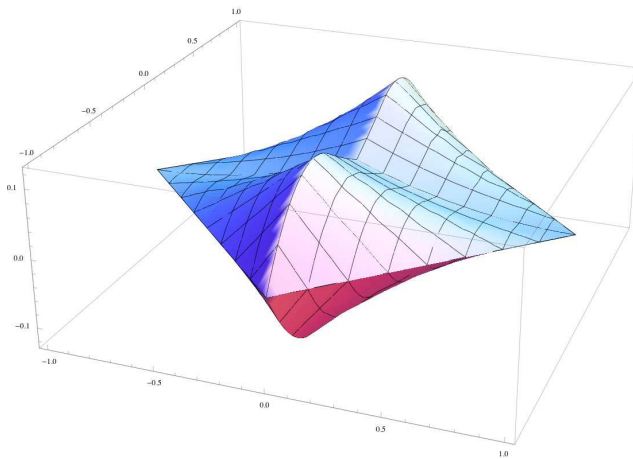}
\caption{
The section of the ball in the plane $y=0$; note the cusps where the vertical direction is squashed.
}
\end{center}
\end{figure}

In \cite{Oberwolfach2006} Burago and Margulis had conjectured that on any discrete group $\Gamma$ any two left-invariant word metrics $\rho_1$ and $\rho_2$ satisfying $\frac{\rho_1(\id,\gamma)}{\rho_2(\id,\gamma)} \to 1$, as $\gamma \to \infty$, must be at a bounded distance from each other, namely $|\rho_1(\id,\gamma) -\rho_2(\id,\gamma)|\leq C$ for all $\gamma \in \Gamma$. This is certainly the case in $\Z^d$ and Krat \cite{krat} established it for the Heisenberg group and for word hyperbolic groups. Abels and Margulis proved an analogous result for word metrics on reductive Lie groups. However, the above example shows that it fails in general.\\

It turns out that a much stronger property holds for the word metrics $\rho_1$ and $\rho_2$ defined above on $G(\Z)$. We prove:
\begin{theorem}\label{noquasi} Even though the two Cayley graphs are quasi-isometric and have isometric asymptotic cones, there is no $C>0$ such that $(G(\Z),\rho_1)$ is $(1,C)$-quasi-isometric to $(G(\Z),\rho_2)$.\\
\end{theorem}

Recall that a map $\phi:X \to Y$ between two metric spaces $(X,d_X)$ and $(Y,d_Y)$ is called a $(1,C)$-quasi-isometry if $d_X(a,b) - C \leq d_Y(\phi(a),\phi(b)) \leq d_X(a,b) +C$, for all $a,b \in X$, and every $y \in Y$ is at distance at most $C$ from some element of $\phi(X)$.\\

This theorem is in sharp contrast with what happens in the Abelian case, where it is a simple matter to establish that two word metrics on $\Z^d$ have isometric asymptotic cones if and only if they are $(1,C)$-quasi-isometric for some $C>0$.\\

Let $d_\infty$ be the Pansu limit metric on the asymptotic cone of $(G(\Z),\rho_i)$. It is easy to see that $$d_\infty(\id,(t,g))=|t|+d_3(\id,g),$$ where $d_3$ is the subFinsler metric on $H_3(\R)$ associated to the $\ell^1$ norm $|x|+|y|$ on the horizontal subspace $\R X \oplus \R Y$. See Figure 1 for a picture of the unit ball of $(H_3(\R),d_3)$.\\

Theorem \ref{noquasi} is a simple consequence of the following proposition, which shows the sharpness of the exponent $\frac{1}{2}$ in Theorem \ref{main-theorem} for step-$2$ groups.\\

\begin{proposition}[Sharpness of $\alpha_2=\frac{1}{2}$]\label{quasi} Let $X^1_n=(B_{\rho_1}(\id,n),\frac{1}{n}\rho_1)$ and  $X_\infty:=(B_{d_\infty}(\id,1),d_\infty)$. Then there is $c>0$ such that
$$d_{GH}(X^1_n,X_\infty) > \frac{c}{\sqrt{n}}.$$\\
\end{proposition}

By contrast the convergence for $X_n^2=(B_{\rho_2}(\id,n),\frac{1}{n}\rho_2)$ is in $O(\frac{1}{n})$, hence much faster.\\

The proof of Proposition \ref{quasi} relies on some geometric considerations pertaining to the precise form of geodesics in the asymptotic cone $(G(\R),d_\infty)$ and some of its finer geometry. The key to it of course is the existence of the abnormal geodesic $\{(t;0,0,0)\}_{t\in[0,1]}$ in $G(\R)$. In the next two sections, we give a sketch of the proof.\\

  \begin{figure}\label{palla2b}
  \begin{center}
\includegraphics[scale=.7]{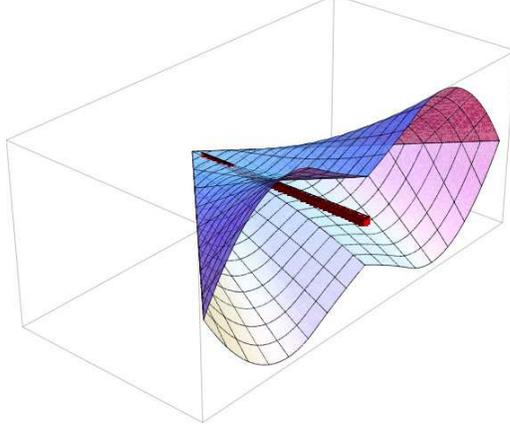}
\caption{
Half of the section of the ball with the abnormal geodesic in red.
}
\end{center}
\end{figure}



\section{Fine geometry of the Heisenberg group equipped with the Pansu metric}\label{fine}
We discuss here the geometry of the Pansu limit metrics on $H_3(\R)$ and $\R \times H_3(\R)$, and state the geometric ingredients needed for Proposition \ref{quasi}.\\

The asymptotic cone of the discrete Heisenberg group $H_3(\Z)$ endowed with standard generators $\{(1,0,0)^{\pm 1}, (0,1,0)^{\pm 1}\}$ is the real Heisenberg group $H_3(\R)$ endowed with the subFinsler metric $d_3$ induced by the $\ell^1$ norm $|x|+|y|$ on the horizontal subspace $\R X\oplus \R Y$ of the Lie algebra. A picture of its unit ball was given in Figure 1. This picture is borrowed from \cite{breuillard}, where we computed the precise form of geodesics in $(H_3(\R),d_3)$.\\

Geodesics in $(H_3(\R),d_3)$ are horizontal paths and can thus be described accurately by their projection to the $(x,y)$-plane, say $(x(t),y(t))$.  There are three kinds of geodesics between $\id$ and a point $g=(x,y,z) \in H_3(\R)$.\\

\begin{enumerate}
\item geodesics of ``staircase type'' where $x(t)$ and $y(t)$ are both monotone (see Figure 4, curve $c$). This happens if and only if $|z|\leq \frac{|xy|}{2},$
\item $3$-sided arcs of square with sides parallel to the $x$-axis and $y$-axis (see Figure 4, curve $a$). This happens if and only if $\frac{|xy|}{2} < |z|\leq \max\{|x|,|y|\}^2 -\frac{|xy|}{2}$,
\item $4$-sided arcs of square with sides parallel to the $x$-axis and $y$-axis (see Figure 4 curve $b$). This happens if and only if $|z|> \max\{|x|,|y|\}^2 -\frac{|xy|}{2}$.
\end{enumerate}

\begin{figure}\label{subFin-geo5}
\begin{center}
\includegraphics[scale=.4]{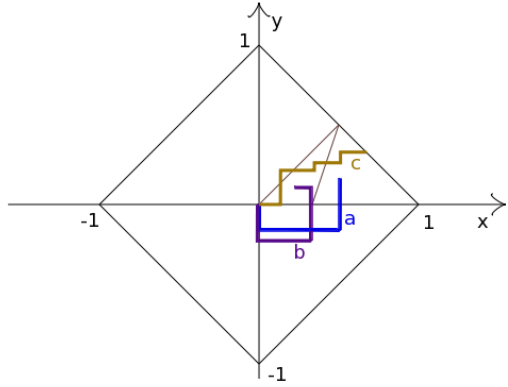}
\caption{
 Geodesics in the Pansu metric $d_3$ on the Heisenberg group $H_3(\R)$.
 There are three kinds of geodesics.
The curve $a$ is the projection of  an example of a geodesic with $3$ sides.
It connects $(0,0)$ to a point in the triangle with vertices $(1,0), (1/3,0), (1/2, 1/2)$.
The curve $b$ is the projection of  an example of a geodesic with $4$ sides.
It connects $(0,0)$ to a point in the triangle with vertices $(0,0), (1/3,0), (1/2, 1/2)$.
The curve $c$ is the projection of   an example of a geodesic of staircase type. It connects $(0,0)$ to a point $(x,1-x)$ with $x\in(0,1)$.
   }
\end{center}
\end{figure}


This classification follows easily from the solution to Dido's isoperimetric problem in the plane equipped with $\ell^1$ norm (see \cite{busemann} and the Appendix to \cite{breuillard}). The uniqueness issue for geodesics is easily addressed: geodesics of staircase type between $\id$ and $g$ are never unique unless $|z|=|xy|/2$. The $3$-sided arcs of square are unique and so are the $4$-sided ones except if $x$ or $y$ is $0$.\\

For more information on the geometry of polygonal subFinsler metrics on the Heisenberg group, we refer the reader to the nice recent preprint by Duchin and Mooney \cite{Duchin-Mooney}.

Accordingly, it is a simple matter to give an exact formula for $d_3$. We obtain:

\begin{enumerate}
\item If $|z|\leq \frac{|xy|}{2}$, then $d_3(\id,(x,y,z))=|x|+|y|$,
\item If $\frac{|xy|}{2} \leq |z|\leq \max\{|x|,|y|\}^2 -\frac{|xy|}{2}$, then $d_3(\id,(x,y,z))= \max\{|x|,|y|\} + \frac{2|z|}{\max\{|x|,|y|\}}$,
\item If $\max\{|x|,|y|\}^2 -\frac{|xy|}{2} \leq |z|$, then $d_3(\id,(x,y,z))=4\sqrt{|z|+\frac{|xy|}{2}}-|x|-|y|$.
\end{enumerate}

The proof of Proposition \ref{quasi} relies on a study of extreme points in the unit balls of $d_3$ and $d_\infty$. A collection of points $g_1,\ldots,g_k$ in the unit ball is said to be a collection of \emph{extreme points} if $d(g_i,g_j)=2$ for every $i \neq j$.\\

It is easy to see that in the unit ball for the $\ell^1$ norm in $\R^k$ there is a unique collection of extreme points of size $2k$, namely the vertices. We prove an analogous characterization of extreme points in $H_3(\R)$ and $\R \times H_3(\R)$:\\

\begin{lemma}[Extreme points]\label{ext-heis} \begin{enumerate}
\item Suppose $g_1,\ldots,g_4$ is a collection of extreme points in the unit ball of $(H_3(\R),d_3)$. Then there are $a,b \in [\frac{1}{2},1]$ such that this collection is $\{(a,1-a,a(1-a)),(1-a,a,-a(1-a)),(-b,-(1-b),b(1-b)),(-(1-b),-b,-b(1-b))\}$, see Figure 5.
\item Suppose $g_1,\ldots,g_6$ is a collection of extreme points in the unit ball of $(\R \times H_3(\R),d_\infty)$. Then this collection is
    $$\{(1;0,0,0),(-1;0,0,0),(0;g_1),(0;g_2),(0;g_3),(0;g_4)\},$$ where $g_1,\ldots,g_4$ are as in item $(i)$.
\end{enumerate}
\end{lemma}

\begin{figure}\label{subFin-geo3}
\begin{center}
\includegraphics[scale=.4]{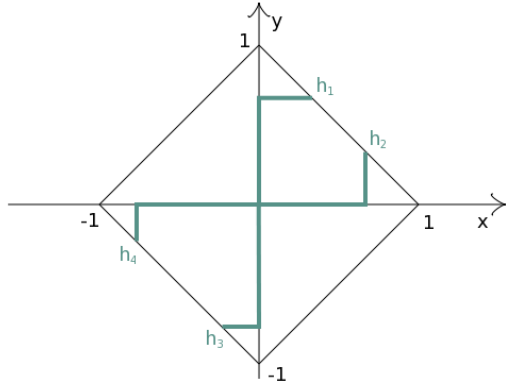}
\caption{A collection of four extreme points in $B_{d_3}(\id,1)$. The extreme points $h_1,\ldots,h_4$ are at distance $2$ from one another and we show the projection of geodesics connecting them.
   }
\end{center}
\end{figure}

This lemma has the following simple consequence for isometries of $X_\infty=B_{d_\infty}(\id,1)$.\\

\begin{lemma}\label{isometry} Any isometry $\phi: X_\infty \to X_\infty$ preserves the pair $\{(1;0,0,0),(-1;0,0,0)\}$.
\end{lemma}

But the above classification of extreme points is crucial to establish the following characterization of almost mid-points of extreme points.\\

\begin{lemma}[Almost mid-points of extreme points]\label{heisen-reduce}
Suppose $h_1, \ldots, h_4$ are $4$ points in $B_{d_3}(\id,1)$ such that $d_3(h_i,h_j) \geq 2- \eps$ for every $i\neq j$. Let $p \in B_{d_3}(\id,1)$ be such that $d_3(\id,p) + d_3(\id,h_i) \leq d_3(p,h_i) + \eps$ for every $i=1,\ldots,4$. Then $d_3(\id,p)=O(\eps)$ as $\eps \to 0$.
\end{lemma}

Finally Lemma \ref{heisen-reduce} is used to establish the following main technical lemma needed in the proof of Proposition \ref{quasi}.\\

\begin{lemma}[Almost extreme points]\label{extreme}
Suppose we are given $5$ points $g_1,\ldots,g_5$ in $B_{d_\infty}(\id,1)$ in $G(\R)=\R \times H_3(\R)$ such that $d_\infty(g_i,g_j) \geq 2 - \eps$ for all $i \neq j$. Let $g=(v;x,y,z) \in B_{d_\infty}(\id,1)$ be such that $d_\infty(g,g_i) \geq 2- \eps$ for every $i=1,\ldots,5$. Then as $\eps \to 0$, either $|v|=O(\eps)$ or $|v-1|=O(\eps)$ or $|v+1|=O(\eps)$.
\end{lemma}

We note that in the unit ball of $H_3(\R)$, $\eps$-extreme points (i.e., points $h_1,...h_4$ as in Lemma \ref{heisen-reduce}) are not necessarily $O(\eps)$ away from genuine extreme points; in general they are only $O(\sqrt{\eps})$ away. The fact that this holds with $O(\eps)$ in the unit ball of $\R \times H_3(\R)$ for the two points close to $(1;0,0,0)$ and $(-1;0,0,0)$ is a manifestation of the presence of the abnormal geodesic $\{(t;0,0,0)\}_{t \in [0,1]}$ and is the heart of the matter here.\\

The proof of Lemmas \ref{ext-heis}, \ref{isometry}, \ref{heisen-reduce} and \ref{extreme} relies on the above classification of geodesics in $(H_3(\R),d_3)$ and the formulas for the distance $d_3$ recalled above.

\section{Proof of Proposition \ref{quasi}}
Set $\eps_n:=d_{GH}(X^1_n,X_\infty)$. By definition of the Gromov-Hausdorff metric, there is a $(1,4\eps_n)$-quasi-isometry $\phi_n:X_\infty \to X^1_n$. Hence there are points $x_n$ and $y_n$ in $X_\infty$ that lie at distance at most $4\eps_n$ from $(1;0,0,\frac{1}{n})$ and $(1;0,0,-\frac{1}{n})$, respectively. But by Theorem \ref{main-theorem-second}, $\rho_1((n;0,0,n),(n;0,0,-n)) \simeq d_\infty((n;0,0,n),(n;0,0,-n))$, i.e., the ratio tends to $1$, as $n \to +\infty$. However, $d_\infty((n;0,0,n),(n;0,0,-n))=d_3(0,0,2n) > c\sqrt{n}$, for some $c>0$. It follows that $d_\infty(x_n,y_n)> c/\sqrt{n} - 8\eps_n$.\\

Therefore we are left to show that $d_\infty(x_n,y_n)=O(\eps_n)$. The sequence $\phi_n$ converges to an isometry $\phi$ of $X_\infty$. By Lemma \ref{isometry} we may assume (up to precomposing all $\phi_n$ by the isometry $v \to -v$) that $\phi$ fixes the point $(1;0,0,0)$. It follows that $x_n$ and $y_n$ converge to $(1;0,0,0)$. Let $\pi$ denote the projection $G(\R) \to \R^3$ modulo the commutator subgroup. Note that $\pi(X_\infty)$ is the $\ell^1$ unit ball in $\R^3$ and has $6$ vertices among which $(1;0,0)$. Considering the $5$ remaining vertices and the map $\pi \circ \phi_n$, we obtain points $g_1,...,g_5 \in X_\infty$ such that $d_\infty(g_i,g_j) \geq 2- \eta_n$ and $d_\infty(g_i,x_n) \geq 2- \eta_n$, where $\eta_n=4\eps_n+\frac{1}{n}$ (observe that $\pi \circ \phi_n$ is an $\eta_n$-submetry). Since $x_n \to (1;0,0,0)$ as $n \to +\infty$, it follows from Lemma \ref{extreme} that $x_n$ is $O(\eta_n)$ close to $(1;0,0,0)$. The same applies to $y_n$ obviously, and hence $d_\infty(x_n,y_n)=O(\eps_n)$ as desired.\\



\section{Asymptotics for the volume of Cayley balls and spheres}
We now pass to the third part of this note and record some applications of our main theorem to volume asymptotics for Cayley balls and Cayley spheres in finitely generated nilpotent groups.\\

In the asymptotic cone, the metric $d_\infty$ scales nicely under the dilation automorphisms, see $(\ref{scaling})$. Hence the volume of balls obeys the law $vol(B_{d_\infty}(\id,t))=Ct^d$, where $C=vol(B_{d_\infty}(\id,1))$. Here $d$ is the Hausdorff dimension of $(G_\infty,d_\infty)$. It is an integer, given by the Bass-Guivarc'h formula \cite{guivarch,bass}:
$$d=\sum_{k \geq 1} kd_k,$$
where $d_k=\dim \g^{(k)}/g^{(k+1)}$. Combined with Theorem \ref{main-theorem} this gives:
\begin{corollary}[Volume asymptotics for balls]\label{volume-theorem} Let $\Gamma$ be a nilpotent group generated by a finite set $S$ with $S=S^{-1}$ and $1 \in S$. Let $B_S(n)=S^n$ be the ball of radius $n$ centered at $\id$ for the word metric $\rho_S$ induced by $S$. Let $r$ be the nilpotency class of $\Gamma$. Then there is $\beta_r>0$ such that 
$$ |B_S(n)| = c_S n^d + O_S(n^{d-\beta_r}), \textnormal{ as } n\to\infty $$
and one can take $\beta_r=1$ if $r\leq 2$ and $\beta_r=\frac{2}{3r}$ if $r>2$.
\end{corollary}

When $\Gamma$ is torsion-free, the constant $c_S$ above is the volume of the unit ball of the asymptotic cone $B_{d_\infty}(\id,1)$ endowed with the Pansu limit metric, where the Haar measure is normalized so that in the Abelianization $\pi(\Gamma)$ has co-volume one in $\pi(G)=\pi(G_\infty)$. Here $\pi: G\to G/[G,G]$.\\

We recall that the asymptotics without error term was proved by Pansu in \cite{pansu} and that the case $r\leq 2$ is a result of Stoll \cite{stoll}. We believe that our error term for $r>2$ is not sharp and that the following holds:
\begin{conjecture}\label{volume-conjecture}
We have $|B_S(n)|=c_S n^d + O_S(n^{d-1})$ for all finitely generated nilpotent groups.
\end{conjecture}

The error term in the volume asymptotics for balls $B_S(n)$ in the Cayley graph of $\Gamma$ is of course related to the volume of spheres $S_S(n)=B_S(n) \setminus B_S(n-1)$. Clearly, if one has the asymptotics $|B_S(n)|=c_S n^d + O(n^{d-\alpha})$ for some $\alpha \leq 1$, then one also have $|S_S(n)|=O(n^{d- \alpha})$. However, the knowledge of an upper bound on the size of the spheres does not seem to give any information on the error terms in the volume of balls.\\

\begin{corollary}[Volume of spheres]\label{volume-spheres} There are constants $C_1,C_2$ depending on $S$ such that,   for all $ n \in \N$ we have
$$C_1 n^{d-1} \leq |S_S(n)| \leq C_2 n^{d-\beta_r},$$
where $\beta_r$ is as in Theorem \ref{volume-theorem}.
\end{corollary}

The upper bound follows immediately from Corollary \ref{volume-theorem}, while the lower bound is a consequence of the following general fact: if $\Gamma$ is any finitely generated group with word metric $\rho_S$, then
$$|B_S(n)| \leq 2n |S_S(n)|.$$

Corollary \ref{volume-spheres} improves on earlier results of Colding and Minicozzi  \cite[Lemma 3.3.]{colding-minicozzi} (also rediscovered by Tessera in \cite{tessera}) bounding from above the volume of spheres in doubling metric spaces in terms of the doubling constant only. Pansu's theorem (i.e., $|B_S(n)| \simeq c_S n^d$) implies that nilpotent groups with word metric $\rho_S$ are doubling metric spaces with doubling constant $\leq (1+\eps)2^d$ for all balls of radius $\geq r(S,\eps)$, and in this case their argument gives an upper bound of the form:
$$ \frac{|S_S(n)|}{|B_S(n)|} = O(n^{-K^{-d}}),$$
for any $K>4$, which is not as good as our Corollary \ref{volume-spheres}, $4^{-d}<\beta_r$ in general.

\section{Concluding remarks}
Observe that Conjecture \ref{stoll-conj} reduces Conjecture \ref{volume-conjecture} to the computation of the asymptotics of the volume of large balls for the Stoll metric $d_S$. Since, unlike $d_\infty$, $d_S$ does not satisfy the nice scaling property \ref{scaling} in general, it is not obvious that the asymptotics of the volume of the balls $B_{d_S}(\id,t)$ has an error term of the form $c_St^d+O(t^{d-1})$.\\

However, one can prove that these balls, when scaled back by the dilation $\delta_{\frac{1}{t}}$, are only $O(\frac{1}{t})$ away from $B_{d_\infty}(\id,1)$ in Hausdorff distance (for a Riemannian metric). Thus the volume asymptotics for $B_{d_S}(\id,t)$ would follow from the following conjectural statement about the unit ball for the Pansu limit metric $d_\infty$ on the asymptotic cone $G_\infty$.\\

\begin{conjecture}[Regularity of subFinsler spheres in Carnot groups]\label{notfractal}
The unit sphere of the Pansu metric $d_\infty$ is rectifiable with respect to any Riemannian distance on $G_\infty$.
In particular, if the group $G_\infty$ has topological dimension $n$, the sphere has  finite $n-1$-dimensional Lebesgue measure.\\
\end{conjecture}

As it turns out, abnormal geodesics are also behind Conjecture $\ref{notfractal}$ above, in fact they are the reason why this conjecture is not obvious and hence neither is the $O(t^{d-1})$ error term in the volume asymptotics of $t$-balls for subFinsler metrics. Indeed, if there were no abnormal geodesics, the distance function $g \mapsto d_\infty(\id,g)$ would be Lipschitz and its level sets (the spheres) would be rectifiable.\\

We recall incidentally that for certain Carnot-Carath\'eodory manifolds, the distance function and the spheres are known to be not subanalytic, see \cite{Bonnard-Chyba-Kupka}.\\

Even if abnormal curves exist in most Carnot groups, they are conjectured to be sparse. According to Montgomery \cite[chapter 10.2]{Montgomery} there ought to be a Sard theorem for the endpoint map, implying in particular that the set of points in $G_\infty$ that can be reached by a singular curve of length at most $1$, say, must be a nowhere dense set of zero Lebesgue measure. This is still an open problem for general Carnot groups. Should the answer be yes, it would then be possible to prove that subFinsler spheres are not fractal objects and that the $n-1$-dimensional Lebesgue measure of subFinsler spheres is finite.\\





\noindent {\bf Acknowledgments.}
E.B. is grateful to the ERC for its support through grant GADA-208091. Both authors would like to thank the MSRI, Berkeley, for perfect working conditions during the Quantitative Geometry program, when part of this research was conducted. E.B. also thanks Fudan University, Shanghai, for its hospitality.

\bibliographystyle{abbrv}
\bibliography{bibfile}











\end{document}